\documentclass{conm-p-l}

\usepackage{amsfonts, amsmath, amssymb, amscd, amsthm}
\usepackage{graphicx}

\newtheorem{thm}{Theorem}[section]
\newtheorem{cor}[thm]{Corollary}
\newtheorem{lem}[thm]{Lemma}
\newtheorem{prop}[thm]{Proposition}

\theoremstyle{definition}
\newtheorem{defn}[thm]{Definition}
\theoremstyle{remark}

\newfont{\eufm}{eufm10}

\begin{document}

\author{D.V. Osin}
\thanks{This work has been supported by the RFFR Grants $\sharp $
02-01-00892, $\sharp $ 03-01-06555.}

\address{Department of Mathematics, 1326 Stevenson Center,
Vanderbilt University, Nashville  TN 37240-0001, USA}
\email{denis.ossine@math.vanderbilt.edu \\
denis.osin@mtu-net.ru}


\title{Weak hyperbolicity and free constructions}

\maketitle

\begin{abstract}
The aim of this note is to show that weak relative hyperbolicity
of a group relative to a subgroup (or relative hyperbolicity in
the sense of Farb) does not imply any natural analogues of some
well-known algebraic properties of ordinary hyperbolic groups. Our
main tools are combination theorems for weakly relatively
hyperbolic groups.
\end{abstract}

\section{Introduction}

Given a group $G$ generated by a set $S$, the {\it Cayley graph}
$\Gamma = \Gamma _S(G)$ of $G$ is an oriented labelled 1--complex
with the vertex set $G$ and the edge set $G\times S$. An edge
$e=(g,s)\in E(\Gamma )$ goes from the vertex $g$ to the vertex
$gs$ and has the label $\phi (e)=s$. The graph $\Gamma $ can be
regarded as a metric space if we endow it with a combinatorial
metric. This means that the length of every edge of $\Gamma $ is
assumed to be equal to 1.

Recall that a geodesic metric space $M$ is hyperbolic, if there
exists $\delta \ge 0$ such that for any geodesic triangle $\Delta
$ in $M$, every side of $\Delta $ is contained in the closed
$\delta $--neighborhood of the union of the other two sides. A
group $G$ is called {\it hyperbolic} if $G$ is generated by a
finite set $X$ and the Cayley graph $\Gamma _X(G)$ is a hyperbolic
metric space \cite{Gro}. This definition is independent of the
choice of the finite generating set $X$. One can generalize the
notion of a hyperbolic group as follows.

\begin{defn}\label{WRH}
Let $G$ be a group, $\mathcal H= \{ H_1, \ldots , H_m\} $ a
collection of subgroups of $G$. A subset $X\subset G$ is a {\it
relative generating set} of $G$ with respect to $\mathcal H$, if
$G$ is generated by $X\cup H_1\ldots \cup H_m$. By the {\it
relative Cayley graph} $\Gamma ^{rel}=\Gamma^{rel}_X(G)$ of $G$
with respect to $\mathcal H$, we mean the Cayley graph of $G$ with
respect to the generating set $ X\cup H_1\cup\ldots \cup H_m$. We
say that $G$ is {\it weakly hyperbolic relative to} $\mathcal H$,
if there exists a finite relative generating set $X$ of $G$ with
respect to $\mathcal H$ and the corresponding relative Cayley
graph is hyperbolic.
\end{defn}

It is straightforward to check  that if $Y$ is another finite
relative generating set of $G$, then the corresponding relative
Cayley graphs $\Gamma ^{rel}_X(G)$ and $\Gamma ^{rel}_Y(G)$ are
quasi--isometric. Since hyperbolicity is preserved under
quasi--isometries, our definition is independent of the choice of
finite relative generating sets of $G$ with respect to $\mathcal
H$. In case the groups $G$, $H_1$, $\ldots $, $H_m$ are finitely
generated in the usual sense, Definition \ref{WRH} is equivalent
to the definition of relative hyperbolicity given by Farb \cite{F}
(see Lemma \ref{EqDef}). Thus Definition \ref{WRH} can be regarded
as a generalization of Farb's one. We use the term 'weakly
relatively hyperbolic group' to distinguish the class of groups
considered in this paper from the class of relatively hyperbolic
groups introduced by Bowditch in \cite{Bow} (the last class is
strictly larger, see \cite{Sz}). It should be noted that Theorems
\ref{main1}, \ref{main2} stated below remain true if we replace
the words 'weakly relatively hyperbolic' with 'relatively
hyperbolic in the sense of Farb' and, in addition, require all
groups and subgroups under consideration to be finitely generated.

It is known that if $G$ is hyperbolic with respect to a collection
of subgroups $\{ H_1, \ldots , H_m\} $ in the sense of Bowditch,
then $G$ inherits some important algebraic and algorithmic
properties of $H_1, \ldots , H_m$ such as finite presentability,
decidability of various algorithmic problems, etc.
\cite{Bum,F,Osin1,Reb}. The main purpose of this note is to show
that this is not so in case of weak relative hyperbolicity. To
construct the corresponding examples we use the following
'combination theorems' for amalgamated products and
HNN--extensions. Other results of this type for hyperbolic and
relatively hyperbolic groups can be found in
\cite{BF,Dah,Git,KM,OM}.

\begin{thm}\label{main1}
Let $H$ be an arbitrary group, $A$ and $B$ two isomorphic
subgroups of $H$. Denote by $G$ the HNN--extension of $H$ with
associated subgroups $A$ and $B$.
\begin{enumerate}
\item $G$ is weakly hyperbolic relative to $H$.

\item If $H$ is weakly hyperbolic relative to $\{ A, B\} $, then
$G$ is weakly hyperbolic relative to $A$.
\end{enumerate}
\end{thm}

\begin{thm} \label{main2}
Let $H$, $K$ be arbitrary groups, $A$ and $B$ isomorphic subgroups
of $H$ and $K$ respectively. Denote by $G$ the amalgamated free
product $H\ast_{A=B} K$.

\begin{enumerate}
\item $G$ is weakly hyperbolic relative to $\{ H, K\} $.

\item If $H$ is weakly hyperbolic relative to $A$ and $K$ is
weakly hyperbolic relative to $B$, then $G$ is weakly hyperbolic
relative to $A$.
\end{enumerate}
\end{thm}

We notice that if $U\le V\le W$ are groups such that $W$ is weakly
hyperbolic relative to $V$ and $V$ is weakly hyperbolic relative
to $U$, then, in general, $W$ is not weakly hyperbolic relative to
$U$. For example, this is so for $U=\mathbb Z$, $V=U\times \mathbb
Z$, $W=V\times \mathbb Z$. Thus the second assertion in each of
these theorems can not be derived from the first ones.

In what follows we call a finitely generated group $G$ {\it
metahyperbolic}, if $G$ is finitely generated and weakly
hyperbolic relative to a hyperbolic subgroup $H\le G$ (or,
equivalently, $G$ is hyperbolic relative to a hyperbolic subgroup
$H\le G$ in the sense of Farb). Such groups are simplest
non--trivial examples of relatively hyperbolic ones. Thus it would
be reasonable to suspect that they are very close to ordinary
hyperbolic groups from the algebraic point of view. However this
is not so as the following examples show.

Recall that any hyperbolic group is finitely presented \cite{Gro}.
Generalizing this fact, the authors of \cite{BC} stated that if a
finitely generated group $G$ is weakly hyperbolic relative to a
finitely generated subgroup $H$, then $G$ is finitely presented
with respect to $H$, which means that a presentation of $G$ can be
obtained from a presentation of $H$ by adding a finite number of
generators and relations. In particular, this would imply that any
metahyperbolic groups is finitely presented in the usual sense.
However, the proof contains a gap which arises from the absence of
the local finiteness of the relative Cayley graph. The following
corollary provides a counterexample.

\begin{cor}
There exists a metahyperbolic group which is not finitely
presented.
\end{cor}

It is well known that any hyperbolic group possess a finite
presentation with Dehn property \cite{Lys}. In particular, the
word problem is decidable for any hyperbolic group in linear time.
Moreover, if $G$ is hyperbolic with respect to a subgroup $H$ and
the word problem is decidable in $H$, then it is decidable in $G$
\cite{F}. (For other algorithmic problems in relatively hyperbolic
groups we refer to \cite{Bum}, \cite{Osin1}, and \cite{Reb}). The
next result shows that this can can not be generalized to the weak
case.

\begin{cor}
There exists a finitely presented metahyperbolic group with
undecidable word problem.
\end{cor}

The last corollary is inspired by the following result. If $G$ is
an infinite hyperbolic group, then $G$ is never simple. Moreover,
if $G$ is not cyclic--by--finite, then it contains uncountably
many normal subgroups \cite{Gro,Ols}. The same is true in case $G$
is relatively hyperbolic in the sense of Bowditch with respect to
an infinite proper subgroup $H$ \cite{Osin2}.

\begin{cor}
There exists a finitely presented infinite metahyperbolic simple
group.
\end{cor}

{\bf Acknowledgements.} The author is grateful to Mike Mihalik and
the referee for useful remarks concerning this paper.

\section{Equivalent definitions of weakly relatively hyperbolic groups}

We begin with various definitions of weak relative hyperbolicity.
Throughout this section we fix a group $G$, a collection of
subgroups $\mathcal H=\{ H_1, \ldots , H_m\} $ of $G$, and a
finite relative generating set $X$ of $G$ with respect to
$\mathcal H$. For a graph $\Xi $, we denote by $V(\Xi )$ and
$E(\Xi )$ the sets of vertices and edges of $\Xi $. If $e$ is an
edge of $\Xi $, we write $e_-$ and $e_+$ for the origin and the
terminus of $e$ respectively.

\begin{defn}
By the {\it left coset graph} $\widetilde \Gamma =\widetilde
\Gamma_X(G)$ of $G$ with respect to $\mathcal H$ we mean the
oriented labelled 1--complex constructed as follows. The vertex
set of $\widetilde \Gamma $ is $V(\widetilde \Gamma )=\{ gH_i,\;
i=1, \ldots , m, g\in G\} $. For two different cosets $fH_i$ and
$gH_j$, there exists an (oriented) edge $e$ going from $fH_i$ to
$gH_j$ if and only if there are elements $a\in fH_i$ and $b\in
gH_j$ such that $b=ax$ for some $x\in X\cup X^{-1}\cup\{ 1\} $.
The triple $(i,j,x)$ is called the {\it label} of $e$. Obviously
$\widetilde \Gamma $ is connected and, in general, not locally
finite.
\end{defn}

The next definition was formulated by Farb \cite{F}.

\begin{defn}
Suppose that the group $G$ is generated by the set $X$ in the
usual (non--relative) sense. We begin with the Cayley graph
$\Gamma _X(G)$ of $G$ and form a new graph as follows: for each
left coset $gH_i$, $i=1, \ldots , m$, of $H_i$ in $G$, add a
vertex $v(gH_i)$ to $\Gamma _X$, and add an edge $e(gh)$ of length
$1/2$ from each element $gh$ of $gH_i$ to the vertex $v(gH_i)$.
The new graph is called the {\it coned--off Cayley graph} of $G$
with respect to $\mathcal H$, and is denoted by $\widehat \Gamma
=\widehat \Gamma _X(G) $.
\end{defn}

We equip the graphs $\widetilde \Gamma $ and $\widehat \Gamma $
with combinatorial metrics. In case $G$ is finitely generated and
$\widehat \Gamma $ is hyperbolic, the group $G$ is called {\it
hyperbolic relative to $\mathcal H$ in the sense of Farb}
\cite{F}.

\begin{defn}\label{QI}
Two metric spaces $M_1, M_2$ are said to be {\it quasi--isometric}
if there exist $\lambda >0$, $c\ge 0$, $\varepsilon \ge 0$, and a
map $\alpha : M_1\to M_2$ such that the following two condition
hold.

\begin{enumerate}
\item[(a)] For any $x,y\in M_1$, we have $$ \frac{1}{\lambda }
dist_{M_1} (x,y)-c \le dist_{M_2} (\alpha (x),\alpha (y))\le
\lambda dist_{M_1} (x,y)+c. $$

\item[(b)] For any $z\in M_2$ there exists $x\in M_1$ such that $$
dist_{M_2}(\alpha (x), z)\le \varepsilon .$$
\end{enumerate}
\end{defn}

The lemma below shows, in particular, that we can regard
Definition \ref{WRH} as a generalization of Farb's one. Recall
that $\Gamma ^{rel} $ denotes the relative Cayley graph of $G$
with respect to $\mathcal H$ defined in the introduction.

\begin{lem}\label{EqDef} The following conditions are equivalent.
\begin{enumerate}
\item[(i)] The graph $\Gamma ^{rel} $ is hyperbolic.

\item[(ii)] The graph $\widetilde \Gamma $ is hyperbolic.
\end{enumerate}
In case $G$ is generated by the set $X$ in the usual
(non--relative) sense, the above conditions are equivalent to
\begin{enumerate}
\item[(iii)] The graph $\widehat \Gamma $ is hyperbolic.
\end{enumerate}
\end{lem}

\begin{proof}
(i)$\Leftrightarrow $(ii).  Recall that hyperbolicity (or the
absence of it) is preserved when we pass from a metric space to a
quasi--isometric one. We define a map $\alpha :V(\Gamma ^{rel})\to
V(\widetilde \Gamma )$ by the rule $\alpha (g)=gH_1$ for any $g\in
G$. Since any graph is quasi--isometric to its vertex set equipped
with the induced metric, it suffices to show that $\alpha $
satisfies conditions (a) and (b) from Definition \ref{QI} for some
$\lambda $, $c$, $\varepsilon $.

Suppose that two vertices $u,v$ are connected by an edge in
$\Gamma ^{rel}$. Then there are only three possibilities. First
assume that there exists $x\in X\cup X^{-1}$ such that $u=vx$.
Clearly $\alpha (u)=uH_1$ and $\alpha(v)=vH_1$ are connected by an
edge in $\widetilde \Gamma $ in this case. Next suppose
$uH_1=vH_1$. Then $\alpha (u)=\alpha (v)$. Finally let $uH_i=vH_i$
for some $i=2, \ldots , m$. Then $\alpha (u)=uH_1$ is connected to
$uH_i=vH_i$ in $\widetilde \Gamma $ by the edge labelled $(1,i,1)$
and $vH_i$ is connected to $\alpha (v)=vH_1$ by the edge labelled
$(i,1,1)$.

Thus in all cases we have $dist_{\widetilde \Gamma }(\alpha
(u),\alpha (v))\le 2$. Obviously this implies
\begin{equation}\label{1}
dist_{\widetilde \Gamma }(\alpha (u),\alpha (v))\le 2
dist_{\Gamma ^{rel}} (u,v)
\end{equation}
for arbitrary $u,v\in \Gamma ^{rel}$. Further, if for some $u,v\in
\Gamma ^{rel}$, $\alpha (u)$ and $\alpha (v)$ are connected by a
path of length $n$ in $\widetilde \Gamma $, then
$v=uh_1x_1h_2x_2\ldots h_nx_n h_{n+1}$, where $h_1, h_{n+1}\in
H_1$, $h_2, \ldots , h_n\in \bigcup\limits_{i=1}^m H_i$, and $x_1,
\ldots ,x_n\in X\cup X^{-1}\cup \{ 1\} $. Therefore, $dist_{\Gamma
^{rel}} (u,v)\le 2n +1$, which yields
\begin{equation}\label{2}
\frac{1}{2} dist_{\Gamma ^{rel}} (u,v)-\frac{1}{2}\le
dist_{\widetilde \Gamma }(\alpha (u),\alpha (v))
\end{equation}
for any $u,v\in \Gamma ^{rel}$. Inequalities (\ref{1}) and
(\ref{2}) together imply the first condition in Definition
\ref{QI} for $\lambda =c=1/2$. It remains to notice that the
second condition holds for $\varepsilon =1$ since any coset $uH_i$
in $\widetilde \Gamma $ is a distance of at most $1$ from
$uH_1=\alpha (u)$.

(ii)$\Leftrightarrow $(iii). Note that the identity map on $G$
induces an isometric embedding $\iota $ of the vertex set
$V(\Gamma ^{rel} )$ of $\Gamma ^{rel} $ to $\widehat \Gamma $ and
$\widehat \Gamma $ belongs to the closed $1$--neighborhood of the
image $\iota (V(\Gamma ^{rel} ))$.
\end{proof}

In the next section we will also use the following result.

\begin{lem}
The group $G$ acts on $\widetilde \Gamma $ by left multiplications
isometrically with a finite number of orbits of edges.
\end{lem}

\begin{proof}
The fact that the action of $G$ is isometric is obvious. Let us
prove that the number of orbits of edges is finite. Note that $G$
can act on $\widetilde \Gamma $ with inversions, so we can not
speak about the quotient of $\widetilde \Gamma $ with respect to
the action.

Clearly it suffices to show that any two edges having equal labels
belong to the same orbit. Since the number of different labels is
finite, this will imply the statement of the lemma. Let $f_1H_i$,
$f_2H_i$, $g_1H_j$, $g_2H_j$ be cosets and $a_1\in f_1H_i$,
$a_2\in f_2H_i$, $b_1\in g_1H_j$, $b_2\in g_2H_j$ elements of $G$
such that $b_1=a_1x$ and $b_2=a_2x$ for some $x\in X\cup
X^{-1}\cup \{ 1\} $. We have to show that there exists $w\in G$
which takes the pair $(f_1H_i, g_1H_j)$ to $(f_2H_i,g_2H_j)$. Let
$w$ be the element $a_2a_1^{-1}$. Obviously $w=f_2hf_1^{-1}$ for
some $h\in H_i$. Thus $w$ takes $f_1H_i$ to $f_2H_i$. Further,
note that
$$w=a_2a_1^{-1}=b_2xx^{-1}b_1^{-1}=b_2b_1^{-1}=g_2kg_1^{-1} $$ for
some $k\in H_j$. Therefore, $wg_1H_j=g_2kH_j=g_2H_j$. This
completes the proof.
\end{proof}

\section{Proofs of the main results}

We start with auxiliary lemmas. The following is a version of
Svar\v{c}--Milnor Lemma for non--proper actions of groups on (not
necessary locally finite) graphs.

\begin{lem}\label{Sh-M}
Let $G$ be a group acting on graphs $\Gamma _1$ and $\Gamma _1$
with finite number of orbits of edges. Suppose that there is a
bijection between the vertex sets $\beta :V(\Gamma _1)\to V(\Gamma
_2)$ such that the diagram
$$
\begin{CD}
V(\Gamma _1) @>\beta >> V(\Gamma _2) \\ @AgAA @AgAA\\ V(\Gamma _1)
@>\beta >> V(\Gamma _2)
\end{CD}
$$
is commutative for any $g\in G$. Then the graphs $\Gamma _1$ and
$\Gamma _2$ are quasi--isometric.
\end{lem}

\begin{proof}
Suppose that two edges $e, f\in E(\Gamma _1)$ belong to the same
orbit, i.e., $ge=f$ for some $g\in G$. Let $p$ be a geodesic paths
connecting $\beta (e_-)$ to $\beta (e_+)$ in $\Gamma _2$.
Obviously $gp$ is a path of the same length as $p$ connecting
$\beta (f_-)$ to $\beta (f_+)$. Thus we have $dist_{\Gamma _2}
(\beta (e_-),\beta (e_+))= dist_{\Gamma _2} (\beta (f_-),\beta
(f_+)).$ Since the number of orbits of edges is finite, there
exists the maximum
$$
M=\max\limits_{e\in E(\Gamma _1)} dist_{\Gamma _2} (\beta
(e_-),\beta (e_+)).
$$
Thus for any $u,v\in V(\Gamma _1)$, we have
$$
dist _{\Gamma _2} (\beta (u), \beta (v))\le M\, dist _{\Gamma
_1}(u,v).
$$
The converse inequality can be obtained in the same way.
Therefore, $\beta $ defines a quasi--isometry between vertex sets
of $\Gamma _1$ and $\Gamma _2$, which yields the assertion of the
lemma.
\end{proof}

Let $\Sigma $ be a graph. For a cycle $c$ in $\Sigma $, we denote
by $[c]$ its homology class in $H_2(\Sigma , \mathbb Z)$. By
$l(c)$ and $d(c)$ we denote the length and the diameter of $c$
respectively. The next proposition is a homological variant of the
characterization of hyperbolic graphs by linear isoperimetric
inequality (see \cite{BowCC,BrH}).

\begin{prop}\label{IP}
For any graph $\Sigma $ the following conditions are equivalent.
\begin{enumerate}
\item[(i)] $\Sigma $ is hyperbolic.

\item[(ii)] There are some positive constants $M$, $L$ such that
if $c$ is a cycle in $\Sigma $, then there exist cycles $c_1,
\ldots , c_k$ in $\Sigma $ with $d(c_i)\le M$ for all $i=1,
\ldots, k$ such that
\begin{equation}\label{c}
[c]=[c_1]+\ldots +[c_k]
\end{equation}
and $k\le Ll(c)$.
\end{enumerate}
\end{prop}

\begin{lem}\label{Ret}
Let $G$ be a group, $R$ a retract of $G$. Suppose that $G$ is
weakly hyperbolic relative to a collection of subgroups $\mathcal
A=\{ A_1,\ldots , A_m\} $ and $A_i\le R$ for all $i=i, \ldots ,m$.
Then $R$ is weakly hyperbolic relative to $\mathcal A$.
\end{lem}

\begin{proof}
Let $X$ be a finite relative generating set of $G$ with respect to
$\mathcal A$, $Y$ the image of $X$ under the retraction $G\to R$.
Then the relative Cayley graph $\Gamma^{rel} _Y(R)$ of $R$ with
respect $\mathcal A$ is a retract of the relative Cayley graph
$\Gamma ^{rel}_X(G)$ of $G$ with respect to $\mathcal A$. Since
hyperbolicity is preserved under retractions, the lemma follows.
\end{proof}

the proofs of the following two lemmas are straightforward and we
left them to the reader.

\begin{lem}\label{AB}
Suppose that a group $G$ is weakly hyperbolic relative to $\{ A,
B, C_1, \ldots , C_m\} $, where $A$ and $B$ are conjugate
subgroups of $G$. Then $G$ is weakly hyperbolic relative to $\{ A,
C_1, \ldots , C_m\} $.
\end{lem}

\begin{lem} Suppose that a group $G$ is weakly hyperbolic relative to
$\{ A, C_1, \ldots , C_m\} $ and $A$ is conjugate to a subgroup
$B$ of $G$. Then $G$ is weakly hyperbolic relative to $\{ B, C_1,
\ldots , C_m\} $.
\end{lem}

Now we are ready to prove Theorems \ref{main1} and \ref{main2}.

\begin{proof}[Proof of Theorem \ref{main1}.]
The group $G$ acts on the Bass--Serre tree $T$. The vertex set
$V(T)$ of $T$ is the set of the left cosets $\{ gH, \; g\in G,\;
i=1, \ldots , m\} $ and $G$ acts on $V(T)$ by left multiplication
\cite{Trees}. Therefore, by Lemma \ref{Sh-M}, the left coset graph
$\widetilde \Gamma $ of $G$ with respect to $H$ is
quasi--isometric to $T$. Thus $G$ is weakly hyperbolic relative to
$H$ by Lemma \ref{EqDef}.

Let us prove the second assertion of the theorem. By Lemma
\ref{AB}, it suffices to show that $G$ is weakly hyperbolic
relative to $\{ A,B\} $. Let
$$G=\langle H,t\; |\; t^{-1}at=\mu (a)\rangle ,$$ where $\mu :A\to
B$ is the isomorphism. By our assumption $H$ is generated by a
finite set $Y$ relative to $\{ A, B\} $. We put $X=\{ t\}\cup Y $.
Obviously $G$ is generated by the finite set $X$ relative to $A$.
Let us consider the relative Cayley graphs $\Gamma ^{rel}_X(G)$ of
$G$ with respect to $\{ A,B\} $ and $\Gamma ^{rel}_Y (H)$ of $H$
with respect to $\{ A,B\} $. We can think of $\Gamma ^{rel}_Y (H)$
as a subgraph of $\Gamma ^{rel}_X(G)$. Let $M$, $L$ be positive
constants such that $\Gamma ^{rel}_Y (H)$ satisfies the second
condition of Proposition \ref{QI}. Without loss of generality, we
may assume $M\ge 4$. For words $U,V$ in the alphabet $X$ we denote
by $\| U\| $ the length of $U$ and write $U\equiv V$ to express
letter--for--letter equality of $U$ and $V$.

Consider a cycle $c$ in $\Gamma ^{rel}_X(G)$. We are going to
check the second condition of the Proposition \ref{IP}. Let $W$ be
the label of $c$, which is a word in the alphabet $X=\{ t\} \cup
Y$. Suppose that $t^{\pm 1}$ appears in $W$ $n$ times. We want to
show that $c$ admits a decomposition of type (\ref{c}) with
$$k\le Ll(c)+n\le (L+1)l(c)$$ terms.

\begin{figure}

\unitlength 1mm 
\begin{picture}(90.38,40)(15,56)
\put(60,74.75){\oval(52.5,31)[]} \put(35.68,61.54){\circle*{.94}}
\put(57.97,90.08){\circle*{.94}} \put(67.49,90.23){\circle*{.94}}
\put(86.22,80.72){\circle*{.94}} \put(86.07,71.8){\circle*{.94}}

\put(61.99,59.2){\vector(-1,0){.07}}
\put(33.72,73.5){\vector(0,1){.07}}
\put(63.92,90.23){\vector(1,0){.07}}
\put(86.30,75.26){\vector(0,-1){.07}}

\put(31.66,73.43){\makebox(0,0)[cc]{$q_1$}}
\put(59.46,57.08){\makebox(0,0)[cc]{$q_2$}}
\put(62.28,93.8){\makebox(0,0)[cc]{$y_1$}}
\put(90.38,76.56){\makebox(0,0)[cc]{$y_2$}}
\qbezier(57.83,90.08)(46.08,74.77)(86.07,71.95)
\qbezier(67.34,90.38)(59.76,80.27)(86.07,80.87)

\put(84.49,87.85){\vector(1,-1){.07}}
\put(61.69,76.54){\vector(2,-1){.07}}
\put(69.42,82.45){\vector(2,-1){.07}}
\put(68.38,79.5){\makebox(0,0)[cc]{$f$}}
\put(59.46,74.62){\makebox(0,0)[cc]{$e$}}
\put(86.51,90.23){\makebox(0,0)[cc]{$v$}}
\end{picture}

\caption{The decomposition of the cycle $c$ in the proof of the
second assertion of  Theorem \ref{main1}.}
\end{figure}
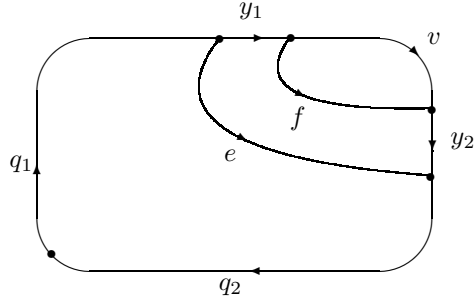

If $n=0$, this is trivial since $W$ represents $1$ in $H$. Further
suppose that $n>0$.  By the Britton Lemma on HNN--extensions (see
\cite{LS}), this means that $W$ has a subword of type $t^{-1}Vt$
of $W$, where $V$ represents an element $a\in A$, or a subword
$tUt^{-1}$, where $U$ represents an element of $B$. We consider
the first case, the second one is analogous. Let $W\equiv
W_1t^{-1}VtW_2$ and let $c=q_1y_1vy_2q_2$ be the corresponding
decomposition of $c$, where $q_1$, $y_1$, $r$, $y_2$ and $q_2$
have labels $W_1$, $t^{-1}$, $R$, $t$, and $W_2$ respectively.
Note that $t^{-1}Vt$ represents some element $b$ of $B$. Thus we
have
\begin{equation}\label{t-string}
[c]=[q_1eq_2]+[e^{-1}y_1fy_2]+[f^{-1}v],
\end{equation}
where $e$ and $f$ are edges of $\Gamma ^{rel}_X(G)$ having labels
$a$ and $b$ respectively (see Figure 1). Note that
$l(q_1eq_2)+l(b^{-1}v)=l(c)$. By the inductive assumption,
$[q_1eq_2]$ and $[f^{-1}v]$ admit decompositions of type (\ref{c})
with at most $Ll(q_1eq_2)+n_1$ and at most $Ll(f^{-1}v)+n_2$ terms
respectively, where $n_1+n_2+2=n$. Together with (\ref{t-string})
this gives a decomposition of type (\ref{c}) for $[c]$ with at
most
$$
k\le Ll(q_1a)+n_1+1+Ll(f^{-1}r)+n_2< Ll(c)+n
$$
terms. Thus $G$ is weakly hyperbolic relative to $\{ A, B\} $ by
Proposition \ref{QI}.
\end{proof}

\begin{proof}[Proof of Theorem \ref{main2}.]
The proof of the first assertion is analogous to that of the first
assertion of Theorem \ref{main1}. Further, recall that the
amalgamated product $G=H\ast _{A=B} K$ is a retract of the
HNN--extension of $H\ast K$ with associated subgroups $A$ and $B$
\cite{LS}. Obviously $H\ast K$ is weakly hyperbolic relative to
$\{ A, B\}$. Thus Theorem \ref{main2} follows from the previous
one.
\end{proof}

\begin{proof}[Proof of Corollary 1.2.]
Let $F$ be a non--cyclic finitely generated free group. Then the
HNN--extension
$$G=\langle F, t\; | \; t^{-1}ft=f,\; f\in [F,F]\rangle $$
is not finitely presented since $[F,F]$ is not finitely generated.
This follows, for example, from the exactness of the
Mayer--Vietoris sequence
$$
\ldots \to H_2(G,\mathbb Z)\to H_1([F,F],\mathbb Z)\to H_1 (F,
\mathbb Z)\to \ldots
$$
or can be proved directly by using the normal form theorem for
HNN--extensions. On the other hand, $G$ is hyperbolic relative to
$F$ by Theorem \ref{main1}.
\end{proof}

\begin{proof} [Proof of Corollary 1.3.]
Recall a result of Rips \cite{Rips}. For any finitely presented
group $Q$ there exists a short exact sequence $$ 1\to N\to H\to
Q\to 1, $$ where $H$ is a finitely generated hyperbolic group and
$N$ is a normal subgroup of $H$ generated by $2$ elements $a,b$.
Let $Q$ be a group with undecidable word problem. Then obviously
the membership problem for $N$ (that is, given an element $h\in
H$, to decide whether $h\in N$) is undecidable. Consider the
HNN--extension
$$G=\langle H,t\;| \; t^{-1}at=a,\; t^{-1}bt=b\rangle ,$$
which is a finitely presented metahyperbolic group. Notice that
$[t,h]=1$ for $h\in H$ if and only if $h\in N$. Therefore, the
word problem is undecidable in $G$.
\end{proof}

\begin{proof} [Proof of Corollary 1.4.]
In \cite{BM}, Burger and Mozes showed that there exists an
infinite simple group $G$ which is an amalgamated product of two
finitely generated free groups along finitely generated subgroups.
Evidently any such a group $G$ is metahyperbolic by Theorem
\ref{main2} since the free group is weakly hyperbolic with respect
to any finitely generated subgroup. (The fact is almost trivial
for free groups; in more general settings it can be found in
\cite{Ger}.)
\end{proof}

\end{document}